\documentclass{aims}
\usepackage{txfonts}
\def\typeofarticle{Type of article}
\def\currentvolume{8}
\def\currentissue{x}
\def\currentyear{2024}

\def\ppages{xxx--xxx}
\def\DOI{10.3934/math.2024xxx}
\def\Received{January 2024}
\def\Revised{February 2024}
\def\Accepted{March 2024}
\def\Published{March 2024}

\numberwithin{equation}{section}

\makeatletter
\renewcommand{\@biblabel}[1]{#1\hfill \hspace{-0.2cm}}
\makeatother

\usepackage{amsmath,amssymb,amsthm}
\usepackage[symbol]{footmisc}
\usepackage{graphicx,hyperref,booktabs,pgf,tikz}
\usepackage{amsrefs}

\usepackage{pgfplots}
\pgfplotsset{compat=1.18} 

\definecolor{graphiccolor}{RGB}{97,129,180}

\theoremstyle{plain}
\newtheorem{thm}{Theorem}

\newtheorem{problem}{Problem}

\theoremstyle{definition}

\newtheorem{example}{Example}

\newtheorem{remark}{Remark}
\theoremstyle{remark}

\begin{document}

\title{The numerical solution of the Dirichlet generalized and
  classical harmonic problems for irregular $n$-sided pyramidal
  domains by the method of probabilistic solutions}

\author{
  Mamuli Zakradze\affil{1}
  Zaza Tabagari\affil{1}
  Nana Koblishvili\affil{1}
  Tinatin Davitashvili\affil{2}
  José-María Sánchez-Sáez\affil{3}
  Francisco Criado-Aldeanueva\affil{4,}\corrauth
}

\shortauthors{
  Zakradze, Tabagari, Koblishvili, Davitashvili, Sánchez-Sáez, Criado-Aldeanueva
}

\address{ 
  \addr{\affilnum{1}}{Department of Computational Methods,
    N. Muskhelishvili Institute of Computational Mathematics of the
    Georgian Technical University, Tbilisi, Georgia.  }
  \addr{\affilnum{2}}{Faculty of Exact and Natural Sciences,
    I. Javakhishvili Tbilisi State University, Tbilisi, Georgia.}
  \addr{\affilnum{3}}{Department of Didactics of Mathematics, Faculty
    of Education, University of Malaga, 29071 Malaga, Spain}
  \addr{\affilnum{4}}{Department of Applied Physics II, Polytechnic
    School, University of Malaga, 29071 Malaga, Spain}
}

\corraddr{Email: fcriado@uma.es}

\begin{abstract}
  This paper describes the application of the method of probabilistic
  solutions (MPS) to numerically solve the Dirichlet generalized and
  classical harmonic problems for irregular $n$-sided pyramidal
  domains. Here, ``generalized'' means that the boundary function has
  a finite number of first-kind discontinuity curves, with the pyramid
  edges acting as these curves. The pyramid's base is a convex
  polygon, and its vertex projection lies within the base. The
  proposed algorithm for solving boundary problems numerically
  includes the following steps: a) applying MPS, which relies on
  computer modeling of the Wiener process; b) determining the
  intersection point between the simulated Wiener process path and the
  pyramid surface; c) developing a code for numerical implementation
  and verifying the accuracy of the results; d) calculating the
  desired function's value at any chosen point. Two examples are
  provided for illustration, and the results of the numerical
  experiments are presented and discussed.
\end{abstract}


\maketitle

\section{Introduction}\label{sec:1}

In practical stationary problems, such as determining the temperature
of a thermal field or the potential of an electric field, it is
sometimes necessary to consider the Dirichlet generalized harmonic
problem (see e.g., \cite{5,3,1,4,2}).

Typically the methods(see e.g., \cite{1,6}) used for obtaining
numerical solutions to classical boundary-value problems are: a) less
suitable or b) ineffective and sometimes useless for solving
generalized boundary problems. The approximate solution of the
generalized problem has low accuracy in the first case due to very
slow convergence of the approximation process near discontinuity
curves (see e.g., \cite{5,3,1,4,2}). The similar
issue occurs in the second case when solving the three-dimensional
Dirichlet generalized harmonic problem using the MFS.

Researchers have attempted preliminary ``improvements'' to address the
boundary value problem at hand. For Dirichlet plane generalized
harmonic problems, the following methods have been developed: I) The
method to reduce the Dirichlet generalized harmonic problem to a
classical problem (see e.g., \cite{7,8}); II) The method of conformal
mapping (see e.g., \cite{9}); and III) The method of probabilistic
solution (see e.g., \cite{11,10}). For
three-dimensional Dirichlet generalized harmonic problems, only the
third method is applicable.
 
Researchers encounter more significant challenges with
three-dimensional Dirichlet generalized harmonic
problems. Specifically, there is no universal approach that can be
applicable to a broad class of spatial domains.
 
The aforementioned references \cite{5,3,1,4,2} addresses relatively
simple problems. Primarily, heuristic methods, method of separation of
variables, particular solutions are used for their solution, resulting
in low accuracy. Heuristic methods can sometimes provide incorrect
solutions, being required verification to ensure all problem
conditions are met (see e.g., \cite{1}). Consequently, the development
of effective computational schemes with high accuracy for the
numerical solution of three-dimensional Dirichlet generalized harmonic
problems, that can be applied to a wide range of spatial domains,
holds both theoretical and practical significance.

It is important to note that in \cite{4}, the discontinuity curves are
disregarded while solving Dirichlet generalized harmonic problems for
a sphere. This oversight, along with the use of classical methods, is
the primary reason for the low accuracy. Thus, for the approximate
solution of three-dimensional Dirichlet generalized harmonic problems,
one should use methods that do not require approximating the
boundary function and are applicable when discontinuity curves
exist. One of such methods is the method of probabilistic solutions
(MPS).

Note that our interest in the numerical solution of both generalized
and classical harmonic Dirichlet problems for pyramidal regions using
the MPS method stems from the following observation: the existing
literature \cite{5,3,1,4,2} does not address either the generalized or
classical problems for pyramidal domains. We believe this omission is
likely due to the geometric complexity of the pyramid.

This observation sparked greater interest in the topic. We began by
developing an appropriate approach for the approximate solution of the
aforementioned problems for a regular $n$-sided pyramid using the MPS
method. Then, we successfully developed the algorithm, which was
applied to solve problems for both a complete pyramid in case of $n=4,
6$, and a regular truncated pyramid in case of $n=5$ (see \cite{12}).

We continued our research in case of irregular $n$-sided pyramids to
solve the same problems using the MPS method. Specifically, we
examined two cases:
\begin{enumerate}
\item An irregular 3-sided pyramid, where the height coincides with a
  lateral edge, lying along the $Ox_3$ axis, with the base on $Ox_{1}
  x_2$ plane, and the angle between side faces containing the height
  is $\pi/2$ (in all cases considered, the pyramid's edges are
  supposed to be the lines of discontinuity);
\item An irregular 4-sided pyramid with a rectangular base located in
  the first quadrant of the $Ox_{1}x_2$ plane, where the vertex
  projection lies at the center of base.
\end{enumerate}
Numerical solutions of these examples using MPS revealed that, unlike
regular pyramids (see \cite{12}), solving generalized and classical
Dirichlet harmonic problems for irregular pyramids requires an
individual approach (see \cite{13}).

In \cite{14}, the problem of finding the numerical solution for both
generalized and classical harmonic Dirichlet problems in a spatially
finite region using the MPS method, particularly at points near the
boundary, and the impact of quantization number $nq$ on solution
accuracy, is studied. The research focused on a specific type of
irregular 3-sided pyramidal region, where all edges are first-order
discontinuity lines (as in the previously discussed examples), and the
height coincides with a lateral edge. Additionally, the angle between
faces containing the height is $\pi/4$. Numerical experiments
demonstrated that achieving an accurate numerical solution near the
boundary requires an optimal selection of the quantization number $nq$
for each point, balancing both solution accuracy and computational
efficiency.

Finally, we developed an appropriate algorithm for solving the
numerical problems using the MPS method, followed by conducting
numerical experiments. The results of this work are presented in this
paper.

The present paper is structured as follows: Section \ref{sec:2}
presents the formulation of the three-dimensional Dirichlet
generalized harmonic problem. Section \ref{sec:3} briefly describes
the MPS and Wiener process simulation. In Section \ref{sec:4} the
algorithm for finding the intersection point between the simulated
Wiener process trajectory and the surface of pyramid is detailed. An
auxiliary classical problem is examined in Section
\ref{sec:5}. Section \ref{sec:6} presents the results of numerical
examples. Finally, Section \ref{sec:7} offers conclusions and
suggestions for future research.

\section{Formulation of the Generalized Problem}\label{sec:2}

Let $D$ be the interior of an irregular $n$-sided pyramid
$P_{n}(h)\equiv P_{n}$ in $R^{3}$ space. By $h$ is denoted its height,
passing through the base of $P_n$. Without loss of generality we
assume that $h$ is aligned with $Ox_{3}$ in the right-handed Cartesian
coordinate system $Ox_{1}x_{2}x_{3}$ and $P_{n}$'s base lies in the
plane $Ox_{1}x_{2}$.

Additionally, suppose the vertices of the base $A_1,A_2,\ldots,A_n$ of
$P_n=MA_1A_2...A_n$ are arranged in a counter-clockwise direction,
with the axis $Ox_{1}$ passing through vertex $A_1$. We will now
formulate the following problem for the pyramid $P_n\equiv \overline
{D}$. Similar to \cite{12}, we now state the following problem for the
irregular pyramid $P_{n}$.

\begin{problem}\label{problem:1}
  Suppose the function $g(y)$, defined on the boundary $S$ of the
  pyramid $P_{n}$, is continuous everywhere except along the edges
  $l_{1},l_{2},\ldots,l_{2n}$, of $P_{n}$, which represent first kind
  discontinuity curves for $g(y)$. The goal is to find a function
  $u(x)\equiv u(x_{1},x_{2},x_{3})\in C^{2}(D)\cap C(\overline
  {D}\backslash \cup_{k=1}^{2n}l_{k})$ that satisfies the following
  conditions:
  \begin{align}
    \Delta u(x)&=0, \;\;\; x\in D, \label{eq:2.1}\\
    u(y)&=g(y), \;\;\; y \in S,\quad y \overline{\in}\; l_k\subset S \quad (k= 1,\ldots,2n),   \label{eq:2.2}\\
    |u(x)|&<c, \;\;\; x\in \overline{D},\label{eq:2.3}
  \end{align}
  where $\Delta=\sum_{i=1}^3\partial^2/\partial x^2_i$ is the Laplace
  operator, $l_{k}$ ($k=1,\ldots,2n$) represent the edges of $P_{n}$,
  and $c$ is a real constant.
\end{problem}

It has been demonstrated (see \cite{16,15}) that the problem
(\ref{eq:2.1}), (\ref{eq:2.2}), (\ref{eq:2.3}) has a solution, which
is unique and continuously depends on the given data. Moreover, there
holds the generalized extremum principle for the generalized solution
$u(x)$:
\begin{equation}
  \min_{x\in S} u(x) <
  \underset{x\in D}{u(x)} <
  \max_{x\in S} u(x),\label{eq:2.4}
\end{equation}
where it is assumed that $x \overline {\in}l_k (k=1,\ldots,2n)$, for
$x\in S$.

Observe (refer to \cite{15}) that the added (\ref{eq:2.3}) condition
of boundedness applies specifically to the neighborhoods surrounding
discontinuity curves of the function $g(y)$ and is crucial for the
extremum principle (\ref{eq:2.4}).

According to (\ref{eq:2.3}), $u(y)$ values are typically undefined
along the $l_k$ curves. For example, if Problem \ref{problem:1} is
concerned with determining the thermal or electric field, then we must take
$u(y)=0$ for $y \in l_k$. In this context, the $l_k$ curves represent
non-conductors (or dielectrics) physically. Otherwise, they would not
be considered as discontinuity curves.

It is evident that, in this approach, the boundary function $g(y)$ has
the following form
\begin{equation}
  g(y)=
  \begin{cases}
    g_1(y),& y\in S_1\\
    g_2(y),& y\in S_2\\
    \vdots &  \vdots\\
    g_{n}(y),& y\in S_{n}\\
    g_{n+1}(y),& y\in S_{n+1}\\
    0,& y\in l_{k}\; (k=1,\ldots,2n)
  \end{cases}
  \label{eq:2.5}
\end{equation}
where: $ S_i$ ($i=1,\ldots,n$) are the lateral faces and $S_{n+1}$ is
the base of $P_n$, respectively, excluding the boundaries. The
functions $g_{i}(y)$, $y\in S_i$ ($i=1,\ldots,n+1$) are continuous on
their respective $S_i$ of $S$. Clearly, $S=(\cup_{i=1}^{n+1}S_i)\cup
(\cup_{k=1}^{2n} l_k)$.

\begin{remark}\label{remark:1}\ 
  \begin{enumerate}
  \item[a)] The case of vacuum within the surface $S$ corresponds to
    the generalized problem with respect to closed shells.
  \item[b)] Not all edges of pyramid must be dielectric in Problem
    \ref{problem:1}. Furthermore, scenarios can be considered where
    faces, apothems, base diagonals, and other elements function as
    dielectrics.
  \end{enumerate}
\end{remark}

\section{Simulation of the Wiener Process and the Method of Probabilistic Solution}\label{sec:3}

This section provides a brief overview of the algorithm for
numerically solving problems of type 1, with a detailed description
available in \cite{17}. The key theorem that enables the application
of the MPS is as follows (see e.g., \cite{16}).

\begin{thm}\label{thm:1}
  Suppose $g(y)$ is continuous (or discontinuous) bounded function on
  $S$ and the finite domain $D\subset \mathbb{R}^3$ is bounded by
  piecewise smooth surface $S$, then the solution of the Dirichlet
  classical (or generalized) boundary value problem for the Laplace
  equation at the fixed point $x\in D$ is given by the following form
  \begin{equation}
    u(x)=E_xg(x(\tau )). \label{eq:3.1}
  \end{equation}
\end{thm}

In equation (\ref{eq:3.1}), it is assumed that the Wiener process
starts at the point $x(t_0)=(x_1(t_0),x_2(t_0),x_3(t_0))\in D$ , where
the value of the desired function is to be determined. $E_x
g(x(\tau))$ represents the mathematical expectation of values of the
boundary function $g(y)$ at the random intersection points between the
Wiener process trajectory and the boundary $S$ and $\tau$ denotes the
random moment at which the Wiener process
$x(t)=(x_1(t),x_2(t),x_3(t))$ first exits the domain $D$.  If $N$, the
number of the random intersection points $y^i=(y_1^i,y_2^i,y_3^i)\in
S\ \ (i=1,2,\ldots ,N)$ is sufficiently large, then by the law of
large numbers, equation (\ref{eq:3.1}) becomes:
\begin{equation}
u(x) \approx  u_N(x)=\frac{1}{N}\sum\limits _{i=1}^N g(y^i) \ \
\label{eq:3.2}
\end{equation}
or $u(x)=\lim u_N(x)$ for $N\rightarrow \infty$, in a probability
sense.

Therefore, in case of the Wiener process the approximate value of the
probabilistic solution to the Problem 1 at a point $x\in D$ is
computed according to the formula (\ref{eq:3.2}).

To simulate the Wiener process, we employ the following recursion
relations (see e.g., \cite{17}):
\begin{equation}
  \begin{gathered}
    x(t_0)=x,\\
    \left.
    \begin{aligned}
      x_1(t_k)&=x_1(t_{k-1})+\gamma_1(t_k)/nq,\\
      x_2(t_k)&=x_2(t_{k-1})+\gamma_2(t_k)/nq,\\
      x_3(t_k)&=x_3(t_{k-1})+\gamma_3(t_k)/nq,
    \end{aligned}
    \right\}
    k=1,2,\ldots,
  \end{gathered}\label{eq:3.3}
\end{equation}
Following these relations, the coordinates of the point $x(t_k) =
(x_1(t_k), x_2(t_k), x_3(t_k))$ are determined. In equation
(\ref{eq:3.3}), $\gamma_1(t_k),\gamma_2(t_k),\gamma_3(t_k)$ are three
independent random numbers with a normal distribution for the $k$-th
step, each with a mean of zero and a variance of one and $nq$
represents the number of quantification where $1/nq =
\sqrt{t_k-t_{k-1}}$. As $nq$ approaches infinity, the discrete process
approximates the continuous Wiener process. In practice, the random
process is simulated at each step of the walk and continues until it
crosses the boundary.

In this paper, we use pseudo-random numbers generated in the MATLAB environment to solve Dirichlet boundary value problems for Laplace's equation.

\section{Finding the Intersection  Point of the  Trajectory of the
  Simulated Wiener Process and the Surface $S$ of the Pyramid
  $P_n$}\label{sec:4}

To determine the intersection points $y^j=(y^j_1,y^j_2,y^j_3)$
($j=1,\ldots,N$) between the Wiener process trajectory and the surface
$S$ (see Section 3), we proceed as follows: First, during the
implementation of the Wiener process, we need to check whether each
current point $x(t_k)$ defined by equation (\ref{eq:3.3}) is inside
$\overline{P_n}$ ($\overline {P_n}=P_n\bigcup S$) or not.

To address this, we have two parameters, $n$ and $h$, as well the
coordinates of vertices: $M=(0,0,h)$, $A_1=(a_1,b_1)$,
$A_2=(a_2,b_2)$,\ldots, $A_n=(a_n,b_n)$ of $P_n$.  Using these
parameters, we need to find: (1) equations of lines passing through
the neighboring vertices of the base; (2) the equations of lateral
faces, and (3) the inclination angles $\alpha_m$, ($m=1,\ldots,n$) of
the lateral faces relative to the base plane of $P_n$.

We will derive the equations of lines that pass through the points
$A_m$ and $A_{m+1}$ ($m=1,\ldots,n$, $A_{n+1}\equiv A_1$) in the
following form
\begin{equation}
  x_2=k_mx_1+c_m.
  \label{eq:4.1}
\end{equation}
Taking into account that $A_m=(a_m,b_m)$, $A_{m+1}=(a_{m+1},b_{m+1})$,
then for definition of constants $k_m$ and $c_m$ in (\ref{eq:4.1}) we
obtain the following algebraic system:
\begin{equation}
  \begin{aligned}
    b_m&=k_ma_m+c_m\\
    b_{m+1}&=k_ma_{m+1}+c_m
  \end{aligned}
  \label{eq:4.2}
\end{equation}
It is easy to see that from (\ref{eq:4.2})
\begin{equation}
  k_m=\frac{b_m-b_{m+1}}{(a_m-a_{m+1})},\qquad
  c_m=\frac{a_mb_{m+1}-b_ma_{m+1}}{a_m-a_{m+1}},
  \label{eq:4.3}
\end{equation}
where if $m=n$, then $a_{n+1}=a_1, b_{n+1}=b_1$.

Without loss of generality, we assume that in (\ref{eq:4.3}),
$a_m-a_{m+1}\neq 0$.  Indeed, if $a_m-a_{m+1}=0$, then the line
$A_mA_{m+1}$ is perpendicular to $Ox_1$.  In order to avoid this
circumstance, it is enough to take such vertex in role of $A_1$ for
which $a_m-a_{m+1}\neq 0$.

Now we find the equations of the planes passing through the points:
$M(0,0,h)$, $A_m$ and $A_{m+1}$, $m=1,\ldots,n$ (or the equations of
the lateral faces of $P_n$).

It is well known that the equation of the plane passing through three
points $P_1(x_1,y_1,z_1)$, $P_2(x_2,y_2,z_2)$, $P_3(x_3,y_3,z_3)$ has
the following form:
\begin{equation}
  \begin{vmatrix}
    x-x_1&y-y_1&z-z_1\\
    x_2-x_1&y_2-y_1&z_2-z_1\\
    x_3-x_1&y_3-y_1&z_3-z_1\\
  \end{vmatrix}=0.
  \label{eq:4.4}
\end{equation}

In our case, in the role of the points: $P_1$, $P_2$, $P_3$ we have
$M(0,0,h)$, $A_m(a_m,b_m,0)$ and $A_{m+1}(a_{m+1},b_{m+1},0)$,
respectively.  On the basis of (\ref{eq:4.4}) the equation of the
plane passing through $M$, $A_m$, $A_{m+1}$ has the following form:
\begin{equation}
  \begin{vmatrix}
    x&y&z-h\\
    a_m&b_m&-h\\
    a_{m+1}&b_{m+1}&-h\\
  \end{vmatrix}
  =0.
  \label{eq:4.5}
\end{equation}

If we write (\ref{eq:4.5}) in system $Ox_1x_2x_3$ and carry out proper
calculations, we will see that the equation has the following form:
\begin{multline}
 h(b_{m+1}-b_m)x_1-h(a_{m+1}-a_m)x_2+(a_mb_{m+1}-a_{m+1}b_m)x_3-\\
 (a_mb_{m+1}-a_{m+1}b_m)h=0.
 \label{eq:4.6}
\end{multline}

In particular, this refers to the equation of the $m$-th lateral face
of $P_n$ when $(x_1,x_2,x_3) \in\overline S_m$.

For the inclination angles $\alpha_m$ ($m=1,\ldots,n$) of the lateral
faces relative to the base plane of $P_n$ we have
$\alpha_m=\arctan(h/d_m)$, where $d_m$ is the distance from the
point $O$ to the $m$-th side of the base of $P_n$ and is given by the
formula: $d_m=\text{abs}(c_m/\sqrt{1+k_m^2}$, ($m=1,\ldots,n$).

With this information about the pyramid $P_n$ one can determine
whether each current point $x(t_k)$, defined by equation
(\ref{eq:3.3}), is inside $\overline{P_n}$ or not. To do this, at each
step of the simulated Wiener process, the angles $\beta_m$
($m=1,\ldots,n$) are calculated, which represent the angles of
inclination of the planes that pass through the points $x(t_k)$,
$A_m$, $A_{m+1}$ relative to the base plane of $P_n$. It is evident
that
\begin{equation*}
  \beta_m=\arctan(x_3(t_k)/\Delta_m),
\end{equation*}
where $\Delta_m$ represents a distance between the point
$(x_1(t_k),x_2(t_k))$ and the line $A_m A_{m+1}$. It is known that
\begin{equation*}
  \Delta_m=\frac{|k_m x_1(t_k)-x_2(t_k)+c_m|}{\sqrt{k_m^2+1}},\qquad
  m=1,\ldots,n,\quad k=1,2,\ldots,
\end{equation*}
where $k_m$ and $c_m$ are defined by (\ref{eq:4.3}).

After calculating the angles $ \beta_m $, they are compared with the
angle $\alpha_m$, ($m=1,\ldots,n$). Specifically:
\begin{itemize}
\item[1)] If $ \beta_m < \alpha_m $ and $ 0<x_3(t_k)<h $ for
  $m=1,\ldots,n$, then the process continues until it intersects the
  boundary $S$
\item[2)] If for $m=p $, $\beta_p=\alpha_p $ and $0<x_3(t_k)<h$ then
  $x(t_k)\in \overline{S_p}$ or $y^j=(y_1^j,y_2^j,y_3^j)=x(t_k)$
\item[3)] If $\beta_p>\alpha_p $ and $0<x_3(t_k)<h $, it indicates
  that the trajectory of the modulated Wiener process at the moment
  $t=t_{k-1}$ has intersected the $p$-th lateral face of $P_n$ or that
  $x(t_{k-1})\in P_n$ while for the moment $t=t_k$, $x(t_k) \overline
  {\in} \overline{P_n}$ . In this case, to approximate the
  intersection point $y^j$, we first determine the parametric equation
  of the line $L$ passing through $x(t_{k-1})$ and $x(t_k)$, which is
  given by the following form:
  \begin{equation}
    \begin{cases}
      x_1=x_1(t_{k-1})+ (x_1(t_k)-x_1(t_{k-1}))\theta,\\
      x_2=x_2(t_{k-1})+ (x_2(t_k)-x_2(t_{k-1}))\theta,\\
      x_3=x_3(t_{k-1})+ (x_3(t_k)-x_3(t_{k-1}))\theta,\\
   \end{cases}
    \label{eq:4.7}
  \end{equation}
  where $(x_1,x_2,x_3)$ is the current point of $L$ and $\theta$ is a
  parameter ($\theta\in\mathbb{R}$).
\end{itemize}

By substituting the expressions for $x_1$, $x_2$ ,$x_3$ from equation
(\ref{eq:4.7}) into equation (\ref{eq:4.6}), we obtain an equation in
terms of $\theta$, which is given by:
\begin{equation}
  A^*\theta=B^*. \label{eq:4.8}
\end{equation}

In (\ref{eq:4.8})
\begin{equation*}
  \left.\begin{aligned}
    A^*=&h(b_{p+1}-b_p)(x_1(t_k)-x_1(t_{k-1}))
    -h(a_{p+1}-a_p)(x_2(t_k)-x_2(t_{k-1}))\\
    &+(a_pb_{p+1}-a_{p+1}b_p)(x_3(t_k)-x_3(t_{k-1}))\\
    B^*=&h(a_pb_{p+1}-a_{p+1}b_p)
    -h(b_{p+1}-b_p)(x_1(t_{k-1})\\
    &+ h(a_{p+1}-a_p)x_2(t_{k-1})
    -(a_pb_{p+1}-b_pa_{p+1})x_3(t_{k-1})
  \end{aligned}\right\}\qquad p=1,\ldots,n.
\end{equation*}

In the given scenario, due to the existence of an intersection point,
$A^*\neq 0$ and $y^j=(x_1(\theta),x_2(\theta),x_3(\theta))$, where
$\theta=B^*/A^*$.

Finally, consider the following cases:
\begin{itemize}
\item[4)] If $x_3(t_k)=0$ and $(x_1(t_k), x_2(t_k))\in S_{n+1}$ then $y^j=
  (x_1(t_k),x_2(t_k),0)$.
\item[5)] If $x_3(t_k)<0$, determine the intersection point
  $y(y_1,y_2,0)$ of the plane $x_3=0$ with the line $L$ that passes
  through the points $x(t_{k-1})$ and $x(t_k)$. Subsequently, if
  $(y_1,y_2,0)\in \overline{S_{n+1}}$ then the intersection point is
  $y^j=(y_1,y_2,0)$.
\end{itemize}

\begin{remark}\label{remark:2}
  It is evident that probability of passing of the path of the
  simulated Wiener process through discontinuous line equals zero, but
  if for any $j$ the intersection point $y^j$ lies on the
  discontinuous line (such case is taken into account in the
  calculation program), then $g(y^j)=0$ is taken in the role of $j$-th
  term of the series in formula (\ref{eq:3.2}).
\end{remark}

\section{An Auxiliary Classical Problem}\label{sec:5}

It is important to note that for three-dimensional case, there are no
exact test solutions available for generalized problems of type
1. Therefore, to verify the accuracy and reliability of the numerical
solution scheme for Problem 1, we use the following approach.

If we set the boundary condition function (\ref{eq:2.5}) in 1 to be
$g_i(y)=1/|y-x^0|$, where $y\in S_i$ ($i= 1,\ldots,n+1$),
$x^0=(x_1^0,x_2^0,x_3^0)\overline \in\overline D$, with $|y-x^0|$
denoting the distance between points $y$ and $x^0$, then the curves
$l_k$ ($k=1,\ldots,2n$) will act as removable discontinuity curves for
the boundary function $g(y)$.  In this case, rather than solving the
generalized Problem \ref{problem:1} we effectively address a classical
Dirichlet harmonic problem.

\begin{problem}\label{problem:2}
  Find a function $u(x)\equiv u(x_1,x_2,x_3) \in C^2(D)\cap
  C(\overline D)$ satisfying the following conditions:
  \begin{equation*}
    \begin{aligned}
      \Delta u(x)&=0, \;\;\; x\in D,\\
      u(y) & =1/|y-x^0|, \qquad \;\;\; y \in S,\quad x^0 \overline \in \overline D.
    \end{aligned}
  \end{equation*}
\end{problem}

We address this problem using the MPS with algorithm designed for
Problem \ref{problem:1}. It is established that Problem
\ref{problem:2} is well posed, meaning that its solution exists, is
unique, and continuously depends on the input data. The exact solution
for Problem \ref{problem:2} is given by
\begin{equation}
  u(x, x^0)=\frac{1}{|x-x^0|},\qquad
  x \in \overline D,\quad
  x^0 \overline \in \overline D.
  \label{eq:5.1}
\end{equation}
In addition to the above, the function $u(x,x^0)$ has the following
properties:
\begin{enumerate}
\item $u(x,x^0)$, as a function of $x$, with $x^0$ fixed, is a
  harmonic function in $\mathbb{R}^3$ everywhere except at the point
  $x^0$;
\item $\lim_{x \to x^0} u(x,x^0)= \infty$;
\item $\lim_{x \to \infty}  u(x,x^0)=0$;
\item In the electrostatic interpretation, $u(x,x^0)$ represents the
  potential at point $x$ in free $\mathbb{R}^3$ space due to a point
  charge of intensity $q=4\pi\,C$ placed at point $x^0$;
\item The function $u(x,x^0)$ is symmetric with respect to its
  arguments in free $\mathbb{R}^3$ space, i.e., $u(x,x^0
  )=u(x^0,x)$. This symmetry is a mathematical representation of the
  principle of reciprocity in physics: a source placed at $x^0$ exerts
  the same influence at point $x$ as a source placed at $x$ does at
  point $x^0$.
\end{enumerate}

Note that, the numerical solution of Dirichlet classical harmonic
problems using the MPS is both intriguing and significant (see
e.g.,\cite{18,20,19}. In the present paper, Problem \ref{problem:2}
serves an auxiliary role, primarily to validate the reliability of the
scheme and the associated program needed for solving Problem
\ref{problem:1}. We first solve Problem \ref{problem:2} and then
compare the results with exact solution of the Problem
\ref{problem:1}, which is solved considering the boundary conditions
(\ref{eq:2.5}).

In this present paper the MPS is applied to two examples. In the
tables in the present paper (see Tables \ref{tab:6.1B} to
\ref{tab:6.2A}), $N$ represents the number of trajectories in the
simulated Wiener process for the given points $x^i=(x_1^i,x_2^i,
x_3^i)\in D$, and $nq$ denotes the quantification number. The tables
below display for problem of type 2, the numerical absolute errors
$\Delta^i $ of $u_N(x)$ as approximated by the MPS, at the points
$x^i\in D$ , for various values of $nq$ and $N$. The results are
presented in scientific notation.

In particular, $\Delta^i=max|u_N(x^i)-u(x^i, x^0)|$,
($i=1,2,\ldots,5$), where $u_N(x^i)$ represents the approximate
solution of Problem \ref{problem:2} at the point $x^i$, as defined by
formula (\ref{eq:3.2}), while the exact solution $u(x^i, x^0)$ of the
test problem is given by equation (\ref{eq:5.1}).  The tables show the
probabilistic solution $u_N(x)$ for Problem \ref{problem:1} at the
points $x^i$, also defined by formula (\ref{eq:3.2}).

\begin{remark}\label{remark:3}
  Problems of types \ref{problem:1} and \ref{problem:2} involving
  ellipsoidal, spherical, cylindrical, conic, prismatic, and regular
  pyramidal domains, as well as axially symmetric finite domains with
  cylindrical holes, external Dirichlet generalized problems for
  spheres, and special types of irregular pyramidal domains, are
  discussed in \cite{15,12,17,22,13,14,21}.
\end{remark}

\section{Numerical Examples}\label{sec:6}

\begin{example}\label{example:6.1}
  In this example the domain $D$ is the interior of an irregular
  3-sided pyramid $P_3(h)$, where $h$ is the height and coordinates of
  the vertices are $M=(0,0,h)$, $A_1=(a_1,b_1)$, $A_2=(a_2,b_2)$,
  $A_3=(a_3,b_3)$.
\end{example}

It is important to note that in these examples, all pyramids are
positioned identically in the coordinate system $Ox_1x_2x_3$ as is
described in Section \ref{sec:2}.

For both types of problems \ref{problem:2} and \ref{problem:1}, we use
$h=2$, $A_1=(3,0)$, $A_2=(0,2)$, $A_3=(-2,-2)$ and $x^0=(0,0,-4)$. In
Problem \ref{problem:1} the boundary function $g(y)\equiv
g(y_1,y_2,y_3)$ is given by
\begin{equation}
  g(y)=
  \begin{cases}
    3,&y\in S_1,\\
    2,&y\in S_2,\\
    1,&y\in S_3,\\
    4,&y\in S_4,\\
    0,&y\in l_k\quad  (k=1,\ldots,6).
  \end{cases}
  \label{eq:6.1}
\end{equation}

In equation (\ref{eq:6.1}), $S_i$ ($i=1,2,3$) and $S_4$ represent the
lateral faces and the base of $P_3$, out of boundaries, respectively,
while $l_k$ ($k=1,\ldots,6$) denote the edges of $P_3$. In this
context, $l_k$ are considered non-conductors (or dielectrics) in a
physical sense.

In all the examples analyzed, we use the scheme described in Section
\ref{sec:4} to determine the intersection points $y^j=(y_1^j, y_2^j,
y_3^j)$ ($j=1,\ldots,N$) where the trajectory of the Wiener process
intersects the surface $S$.

As previously mentioned, we first solve the auxiliary Problem
\ref{problem:2} to validate the accuracy of calculation program for
Problem \ref{problem:1}.

Since, the accuracy of approximate solution $u_N(x)$ of the test
Problem \ref{problem:2} using the MPS at any point $x\in D$ depends
on: complexity of the problem domain; location of point $x$ in the
domain $D$ and values of numbers $nq$ and $N$. Therefore, in
considered examples, during the approximate solution of Problem
\ref{problem:2}, in the first place, for each point $x$, we find the
optimal quantification number $nq$ (by selection) in the sense of
accuracy (see e.g.,\cite{14}).

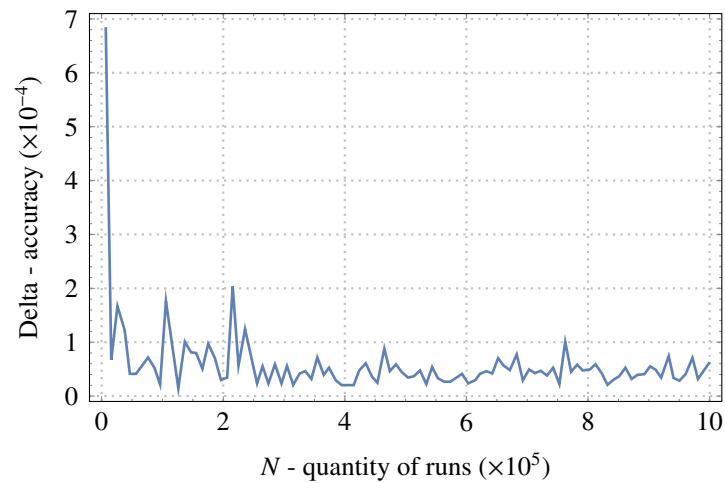
\begin{figure}[htbp]
  \begin{center}
    \caption{Results for the Problem \ref{problem:2} (Example
      \ref{example:6.1}) for starting point $(0, 0, 0.5)$}
    \label{fig:1}
    
    \begin{tikzpicture}
      \begin{axis}[
          xlabel={$N$ - quantity of runs ($\times 10^5$)}, ylabel={Delta - accuracy ($\times 10^{-4}$)},
          grid=major, 
          major grid style={thick,dotted},
          ytick distance=1,
          minor y tick num=4,
          minor x tick num=3,
          enlargelimits=false,
          xmin=-.2,xmax=10.2,
          ymin=-.1,ymax=7.1,
          tick label style={font=\small},
          label style={font=\small},
          major tick length=2pt,
          minor tick length=1pt,
          x=8cm/10,
          y=5cm/7
        ]
        \addplot[color=graphiccolor,line width=1pt] table[x=x, y=y, col sep=comma] {datos.csv};
      \end{axis}
    \end{tikzpicture}
  \end{center}
\end{figure}

Table \ref{tab:6.1B} (see also the Figure \ref{fig:1}) displays the
absolute errors $\Delta^i$ for the numerical solution $u_N(x)$ of the
test Problem \ref{problem:2} at the points $x^i \in D$
($i=1,\ldots,5$).

\begin{table}[htbp]
  \begin{center}
    \caption{Results for the Problem \ref{problem:2} (Example \ref{example:6.1})}\label{tab:6.1B}
    \begin{tabular}{cccccc}\toprule
      $x^i $&$(0,0,0.2)$&$(0,0,0.5)$&$(0,0,1)$&$(0,0,1.5)$&$(0,0,1.8)$\\\midrule
      $N$&$\Delta^1,nq=400$&$\Delta^2,nq=200$&$\Delta^3,nq=300$&$\Delta^4,nq=400$&$\Delta^5,nq=400$ \\\midrule
      $5E+3$&$0.21E-3$&$0.64E-3$&$0.18E-3$&$0.26E-3$&$0.43E-5$ \\
      $1E+4$&$0.17E-3$&$0.11E-3$&$0.11E-4$&$0.53E-4$&$0.51E-4$ \\
      $5E+4$&$0.35E-5$&$0.97E-4$&$0.94E-5$&$0.81E-4$&$0.32E-4$ \\
      $1E+5$&$0.53E-4$&$0.43E-5$&$0.81E-4$&$0.21E-4$&$0.26E-4$ \\
      $5E+5$&$0.44E-4$&$0.10E-5$&$0.39E-4$&$0.42E-4$&$0.19E-4$ \\
      $1E+6$&$0.37E-4$&$0.12E-5$&$0.28E-4$&$0.26E-4$&$0.25E-4$ \\
      \bottomrule
    \end{tabular}
  \end{center}
\end{table}

Based on the numerical results shown in Table \ref{tab:6.1B} we can
confirm that the code for Problem \ref{problem:1} is correct. In Table 1, the interior points are chosen for the numerical experiment. Although any point can be selected, including those from the neighborhood of edges. This issue is discussed in more detail in \cite{23}.

Generally, increasing the values of $nq$ and $N$ enhances the accuracy
of the solution. To achieve this, high performance computing resources are required.

\begin{figure}[htbp]
  \begin{center}
    \caption{ Illustration of single run for the Problem \ref{problem:2} (Example
      \ref{example:6.1})\: ($13657$ steps)}
    \label{fig:2}
    \begin{tikzpicture}

    \end{tikzpicture}
  \end{center}
\end{figure}

For this specific run, 13657 “very small” steps were required from the starting point to the intersection with the boundary of the pyramid. Solving the stated problem requires millions of such runs.

Table \ref{tab:6.1A} presents the numerical solution $u_N(x)$ for
Problem \ref{problem:1} at the same points $x^i$ ($i=1,\ldots,5$)
along with the corresponding values of $nq$ (refer to Table
\ref{tab:6.1B}).  The results demonstrate sufficient accuracy for many
practical applications and align well with the expected physical
behavior.

\begin{table}[htbp]
 \begin{center}
   \caption{Results for Problem \ref{problem:1} (Example \ref{example:6.1})}\label{tab:6.1A}
   \begin{tabular}{cccccc}\toprule
     $x^i $&$(0,0,0.2)$&$(0,0,0.5)$&$(0,0,1)$&$(0,0,1.5)$&$(0,0,1.8)$\\\midrule
     $N$&$u_N(x^1)$&$u_N(x^2)$&$u_N(x^3)$&$u_N(x^4)$&$u_N(x^5)$ \\\midrule
     $5E+3$&$3.42580$&$2.72620$&$1.99640$&$1.82100$&$1.77840$ \\
     $1E+4$&$3.43140$&$2.69770$&$2.00030$&$1.80340$&$1.79950$ \\
     $5E+4$&$3.42110$&$2.69890$&$2.01866$&$1.81518$&$1.79568$ \\
     $1E+5$&$3.42644$&$2.69704$&$2.01578$&$1.81060$&$1.79759$ \\
     $5E+5$&$3.42183$&$2.69627$&$2.01549$&$1.81343$&$1.79700$ \\
     $1E+6$&$3.42202$&$2.69564$&$2.01303$&$1.81354$&$1.79602$ \\
     \bottomrule
   \end{tabular}
 \end{center}
\end{table}

\begin{example}\label{example:6.2}
  In the present example, the domain $D$ represents the interior of an
  irregular 8-sided pyramid $P_8(h)$ with $h=2$. The coordinates of
  the vertices are: $M(0,0,2)$, $A_1=(2,0)$, $A_2=(1.5,1)$,
  $A_3=(0,2)$, $A_4=(-1,2)$, $A_5=(-2,0)$, $A_6=(-1.5,-1)$,
  $A_7=(0,-2)$, $A_8=(1,-1.5)$ and $x^0=(0,0,-4)$.
\end{example}

For Problem \ref{problem:1} the boundary function $g(y)$ is defined as
follows:
\begin{equation}
  g(y)=\begin{cases}
  1,& y\in S_1,\\
  0,& y\in S_2,\\
  2,& y\in S_3,\\
  1,& y\in S_4,\\
  3,& y\in S_5,\\
  0,& y\in S_6,\\
  1,& y\in S_7,\\
  2,& y\in S_8,\\
  4,& y\in S_9,\\
  0,& y\in l_k\  (k= 1,\ldots,16).
  \end{cases}
  \label{eq:6.2}
\end{equation}

In equation (\ref{eq:6.2}), $S_i$ ($i=1,\ldots,8$) and $S_9$ represent
the lateral faces and the base of $P_8$ out of boundaries,
respectively, while $l_k$ ($k=1,\ldots,16$) are the edges of
$P_8$. Additionally, in this case the edges $l_k$ and $S_2$, $S_6$ are
non-conductors.

For Example \ref{example:6.2}, both Problems \ref{problem:2} and
\ref{problem:1} were solved using the same parameters, as in Example
\ref{example:6.1}, specifically the points $x^i$ ($i=1,\ldots,5$) and
numbers $nq$ and $N$.

Table \ref{tab:6.2B} displays the absolute errors $\Delta^i$ of the
numerical solution $u_N(x)$ for the test Problem \ref{problem:2} at
the points $x^i\in D$ ($i=1,\ldots,5$).

\begin{table}[htbp]
  \begin{center}
    \caption{Results for Problem \ref{problem:2} (Example \ref{example:6.2})}\label{tab:6.2B}
    \begin{tabular}{cccccc}\toprule
      $x^i $&$(0,0,0.2)$&$(0,0,0.5)$&$(0,0,1)$&$(0,0,1.5)$&$(0,0,1.8)$ \\\midrule
      $N$&$\Delta_1,nq=400$&$\Delta_2,nq=200$&$\Delta_3,nq=300$&$\Delta_4,nq=400$&$\Delta_5,nq=400$ \\\midrule
      $5E+3$&$0.26E-3$&$0.33E-3$&$0.15E-3$&$0.12E-4$&$0.12E-3$\\
      $1E+4$&$0.24E-3$&$0.24E-3$&$0.31E-3$&$0.81E-4$&$0.22E-4$\\
      $5E+4$&$0.77E-4$&$0.33E-4$&$0.10E-3$&$0.35E-4$&$0.52E-4$\\
      $1E+5$&$0.50E-4$&$0.54E-4$&$0.74E-4$&$0.72E-4$&$0.96E-5$\\
      $5E+5$&$0.49E-4$&$0.21E-4$&$0.40E-4$&$0.57E-4$&$0.29E-4$\\
      $1E+6$&$0.48E-4$&$0.48E-4$&$0.31E-5$&$0.95E-5$&$0.16E-4$\\
      \bottomrule
    \end{tabular}
  \end{center}
\end{table}

Based on the numerical results presented in Table \ref{tab:6.2B}, we
can confirm that the code for the Problem \ref{problem:1} is correct.

Table \ref{tab:6.2A} displays the approximate values of the solution
$u_N(x)$ for the Problem \ref{problem:1} at points $x^i$
($i=1,\ldots,5$).

\begin{table}[htbp]
  \begin{center}
    \caption{Results for Problem \ref{problem:1} (in Example \ref{example:6.2})}\label{tab:6.2A}
    \begin{tabular}{cccccc}\toprule
      $x^i $&$(0,0,0.2)$&$(0,0,0.5)$&$(0,0,1)$&$(0,0,1.5)$&$(0,0,1.8)$\\\midrule
      $N$&$u_N(x^1)$&$u_N(x^2)$&$u_N(x^3)$&$u_N(x^4)$&$u_N(x^5)$\\\midrule
      $5E+3$&$3.41100$&$2.57060$&$1.57640$&$1.13120$&$1.06340$ \\
      $1E+4$&$3.37220$&$2.54410$&$1.59120$&$1.13420$&$1.05030$ \\
      $5E+4$&$3.37734$&$2.56292$&$1.56556$&$1.13714$&$1.06130$ \\
      $1E+5$&$3.38355$&$2.55075$&$1.57076$&$1.13227$&$1.05812$ \\
      $5E+5$&$3.38580$&$2.55173$&$1.56809$&$1.13250$&$1.06473$ \\
      $1E+6$&$3.38650$&$2.55202$&$1.56873$&$1.13109$&$1.06451$ \\
      \bottomrule
    \end{tabular}
  \end{center}
\end{table}

The results are sufficiently accurate for many practical applications and align well with any physical process described by the Laplace equation.

In the present work we solved the problems of type 1 where the
boundary functions $g_i(y)$ ($i=1,\ldots,n+1$) are constants. This
choice was driven by the aim to evaluate how closely the calculation
results match the physical reality. It is clear that solving Problem
\ref{problem:1} under condition (\ref{eq:2.5}) is as easy as solving
Problem \ref{problem:2}. Generally, the Problem \ref{problem:1} can be
solved for configurations of discontinuity curves that allow us to
identify the section of the surface $S$ where the intersection point
lies.
    
The analysis of numerical experiments demonstrates that the proposed
algorithm is both reliable and effective for solving problems of types
1 and 2. In particular, the algorithm is notably straightforward to be
implemented numerically.

\section{Concluding remarks }\label{sec:7}

\begin{enumerate}
\item This study demonstrates that the proposed algorithm is
  exceptionally well-suited to solve approximately the Problems
  \ref{problem:2} and \ref{problem:1}. Notably, this algorithm can
  determine the solution at any point within the domain, a feature
  that distinguishes it from other algorithms discussed in the
  literature.
\item A key advantage of the algorithm is that it does not require an
  approximation of the boundary function, which enhances its
  effectiveness.
\item The algorithm is user-friendly to code, has low computational
  cost, and provides accuracy that meets the requirements of many
  practical problems.
\item In future work, to explore the following areas is planned:
  \begin{itemize}
  \item Applying the proposed algorithm to solve Dirichlet classical
    and generalized harmonic problems in an infinite $\mathbb{R}^3$
    space with a finite number of spherical cavities.
  \item Utilizing the MPS to tackle similar problems within finite
    domains enclosed by multiple closed surfaces.
  \item Extending the MPS application to solve these problems in
    infinite 2D domains with a finite number of circular holes.
  \item Implementing the MPS for solving these problems in regular and
    irregular prisms.
  \end{itemize}
\end{enumerate}

\section*{Author contributions}

Conceptualization, M.Z., Z.T., N.K. and T.D.; methodology, M.Z., Z.T.,
N.K.,T.D., J.M.S. and F.C.-A.; formal analysis,M.Z., Z.T., N.K.,
T.D. and F.C.-A.; investigation,M.Z., Z.T., N.K.,T.D., J.M.S. and
F.C.-A.; writing—original draft,M.Z., Z.T., N.K. and T.D.;
writing—review and editing, J.M.S. and F.C.-A. All authors have read
and agreed to the published version of the manuscript.

\section*{Conflict of interest}

The authors declare no conflict of interest.



\end {document}